\documentclass[12pt,reqno]{amsart}
\usepackage{graphicx, amssymb, hyperref, color,pdfsync}
\usepackage[all,cmtip]{xy}
\numberwithin{equation}{section}
\usepackage{mathrsfs}
\usepackage{tikz}
\usetikzlibrary{matrix,arrows}
\usepackage{dsfont}
\usepackage[USenglish]{babel}
\usepackage{longtable}

\input epsf
\advance\hoffset by -0.9 truecm

\newenvironment{psmallmatrix}
  {\left(\begin{smallmatrix}}
 {\end{smallmatrix}\right)}
\begin{document}
\newcommand{\s}{\vspace{0.2cm}}

\newtheorem{theorem}{Theorem}[section]
\newtheorem{proposition}[theorem]{Proposition}
\newtheorem{prop}[theorem]{Proposition}
\newtheorem{lemma}[theorem]{Lemma}
\newtheorem{lem}[theorem]{Lemma}
\newtheorem{corollary}[theorem]{Corollary}
\newtheorem{assumption}[theorem]{Assumption}
\newtheorem{problem}[theorem]{Problem}
\newtheorem{observation}[theorem]{Observation}
\newtheorem{rema}[theorem]{Remark}
\newtheorem{assertion}[theorem]{Assertion}
\newtheorem{result}[theorem]{Result}
\newtheorem{fact}[theorem]{Fact}
{\theoremstyle{definition}
\newtheorem{definition}{Definition}}
\newtheorem*{theo*}{Theorem}
\newtheorem{example}[theorem]{Example}
\newtheorem{rem}[theorem]{Remark}
\newtheorem{question}[theorem]{Question}
\newtheorem*{rem*}{Remark}
\newtheorem*{corollary*}{Corollary}

\title[Classifying compact Riemann surfaces by number of  symmetries]{Classifying compact Riemann surfaces \\by number of  symmetries}

\author[Reyes-Carocca]{Sebasti\'an Reyes-Carocca} 
\address{Departamento de Matem\'aticas, Facultad de Ciencias, Universidad de Chile, Las Palmeras 3425, \~{N}u\~noa, Santiago, Chile.}
\email{sebastianreyes.c@uchile.cl}

\author[Speziali]{Pietro Speziali}
\address{Instituto de Matem\'atica, Estat\'istica e Computa\c{c}\~ao Cient\'ifica, Universidade Estadual de Campinas, Campinas, SP 13083-859, Brazil.}
\email{speziali@unicamp.br}

\keywords{Compact Riemann surfaces, group actions, automorphisms}
\subjclass[2010]{30F10, 32G15, 14H37, 30F35, 14H30}
\thanks{Partially supported by Fondecyt Regular Grant 1220099 and 1230708, and MATH-AMSUD Grant 22-MATH-03.}

\begin{abstract} In this article we consider compact Riemann surfaces that are uniquely determined by the property of possessing a group of automorphisms of a prescribed order,  strengthening  uniqueness results proved by Nakagawa. More precisely, we deal with the cases in which such an order is $3g$ and $3g+3,$ where $g$ is the genus. We prove that if $g$ is odd (respectively $g$ even and $g \not \equiv 2 \mbox{ mod }3$) then there exists a unique  Riemann surface of genus $g$ with a group of automorphisms of order $3g$ (respectively $3g+3$). A similar conclusion can be derived in terms of orientably-regular hypermaps. In addition, we determine the full automorphism group of such Riemann surfaces  and provide decompositions of their Jacobians.
\end{abstract}

\maketitle
\thispagestyle{empty}
\section{Introduction}
Let $\mathscr{M}_g$ denote the moduli space of  compact Riemann surfaces (or smooth complex projective algebraic curves) of genus $g \geqslant 2.$ It is classically known that if $X \in \mathscr{M}_g$ then $X$ has finitely many automorphisms, and that 
 \begin{equation*}\label{eq:h}
 |\mbox{Aut}(X)| \leqslant 84(g-1).
 \end{equation*}This upper bound is called the Hurwitz bound.
 
Riemann surfaces endowed with many symmetries have been extensively studied in recent years, and up to the present there are still many open problems on this beautiful subject; see, for instance, the expository papers  \cite{BPW22}, \cite{Bu21} and \cite{survey2}.
\s
 
By the Hurwitz bound, it is natural to consider Riemann surfaces endowed with a group of automorphisms whose order is $\mathbb{Q}$-linear in the genus of the surface. More precisely, following \cite{K911}, a sequence $ag+b$, where $a$ and $b$ are rational numbers is called admissible if for infinitely many values of $g$ there is a compact Riemann surface of genus $g$ with a group of automorphisms of order $ag + b$. Once an admissible sequence is given, the problem of classifying and describing such Riemann surfaces arises naturally. 

\s

A classical example is the admissible sequence $8g+8$. It is known that there is a unique, up to isomorphism, compact Riemann surface of genus $g$ with a group of automorphisms of order $8g+8$ if and only if $g \not\equiv 3 \mbox{ mod } 4$ is sufficiently large. This Riemann surface is represented by the so-called Accola-Maclachlan  curve $$y^2=x^{2g+2}-1.$$We refer to the articles \cite{Accola}, \cite{K91} and \cite{Mac} for more details.   This well-known example shows that one of the challenges of studying this kind of problem is not only to find an infinite family of genera for which the uniqueness holds, but also to find the largest set of genera for which such  uniqueness is achieved.

\s

 In this paper, we successfully manage to do so for the admissible sequences $3g$ and $3g+3.$ Such admissible sequences arise from the study of cyclic groups of automorphisms of maximal possible order, as we explain below.

\s

Let $g \geqslant 2$ and $c \geqslant 1$ be fixed integers. We denote by $\mathscr{N}_{g,c}$ the maximum of the orders bounded by $c$ of automorphisms of Riemann surfaces of genus $g$, namely, $$\mathscr{N}_{g,c}=\max_{X \in \mathscr{M}_g}\{\mbox{o}(h): h \in \mbox{Aut}(X) \mbox{ and } \mbox{o}(h) < c\}.$$

A well-known result of Wiman \cite{Wi} (see also \cite{Harvey66}) is the fact that $$\mathscr{N}_{g,\infty}=\max_{X \in \mathscr{M}_g}\{\mbox{o}(h): h \in \mbox{Aut}(X)\} =4g+2$$and that this maximum is attained by the compact Riemann surface  of genus $g$ algebraically represented by the so-called Wiman curve of type I $$W_{g,1}: y^2=x^{2g+1}-1.$$

Later, Nakagawa in \cite{N84} (see also \cite{K97}) considered the uniqueness problem and proved that a Riemann surface of genus $g$ endowed with an automorphism of order $4g+2$ is necessarily isomorphic to $W_{g,1}.$  Soon after, Kulkarni in \cite{K91} took this uniqueness problem further by considering the question as to how far is $W_{g,1}$ from being determined by the order of its automorphism group alone. Concretely, he succeeded in proving  that if $g$ is sufficiently large and $2g+1$ is not of the form $3m$ or $9m$ where all prime divisors of $m$ are congruent to 1 mod 3, then
$$X \in \mathscr{M}_g \mbox{ and } |\mbox{Aut}(X)|=4g+2 \, \iff \, X \cong W_{g,1}.$$

Wiman in \cite{Wi} also showed that the so-called Wiman curve of type II  $$W_{g,2}: y^2 = x(x^{2g} -1)$$represents a compact Riemann surface of genus $g$ admitting an automorphism of order $4g$.  Furthermore, as proved in \cite{K97}, for each $g \geqslant 2$ one has that $$\mathscr{N}_{g,4g+2}=4g,$$ and  if $g \neq 3$ then $W_{g, 2}$ is the only Riemann surface of genus $g$ admitting an automorphism of order $4g$.  If $g=3$, then the Picard curve $y^3 = x^ 4 - 1$ also has an automorphism of order twelve.  We refer to \cite{N84} and \cite{K97} for more details.

\s

In contrast with the case $4g+2,$ there are infinitely many compact Riemann surfaces of genus $g$ with a group of automorphisms of order $4g$. Indeed, following \cite{BCI17}, up to five exceptional values of $g$, all  compact Riemann surfaces of genus $g$ with exactly $4g$ automorphisms form a  complex one-dimensional equisymmetric family (see also \cite{SRC19}).

\s

The computation of $\mathscr{N}_{g,4g}$, namely, the largest order of automorphisms of Riemann surfaces of genus $g$ that are bounded by $4g$,  turns out to be irregular and depends on the residue of $g$ modulo three, as proved by Nakagawa in \cite{N84}. For the sake of explicitness and  to introduce some needed notation, we review these cases here.
\subsection*{1.} Assume $g \equiv 0 \mbox{ mod }3$. Then $\mathscr{N}_{g,4g}=3g+3$ and there is a unique Riemann surface of genus $g \neq 3$ with an automorphism of order $3g+3,$  represented by the curve  $$C_{g,1}:y^3=x^2(x^{g+1}-1).$$In addition, $\mathscr{N}_{g,3g+3}=3g$ and there is a unique Riemann surface of genus $g \neq 6,9,12$ with an automorphism of order $3g,$  represented by the curve  $$C_{g,2}:y^3=x(x^{g}-1).$$
\subsection*{2.}  Assume $g \equiv 1 \mbox{ mod }3$. Then $\mathscr{N}_{g,4g}=3g+3$ and there is a unique Riemann surface of genus $g$ with an automorphism of order $3g+3,$  represented by the curve $$C_{g,3}:y^3=x(x^{g+1}-1).$$In addition, $\mathscr{N}_{g,3g+3}=3g$ and there is a unique Riemann surface of genus $g \neq 4,10$ with an automorphism of order $3g,$  represented by the curve $C_{g,2}.$

\subsection*{3.}  Assume $g \equiv 2 \mbox{ mod }3.$  Then $\mathscr{N}_{g,4g}=3g$ and there is a unique Riemann surface of genus $g \neq 2$ with an automorphism of order $3g,$  represented by the curve $$C_{g,4}:y^3=x^2(x^{g }-1).$$

\s

In this article, we consider compact Riemann surfaces of genus $g$ endowed with (not necessarily cyclic) groups of automorphisms of order $3g$ and $3g+3.$ More precisely, we address the problem of determining the extent to which the before introduced Riemann surfaces are uniquely  determined by the order of their automorphisms groups alone.

\s

Our results can be succinctly summarised as follows.

\s

{\bf a.} We prove that, for each $g \neq 21$ odd, $C_{g,2}$ and $C_{g,4}$ are the only Riemann surfaces of genus $g$ with a group of automorphisms of order $3g.$ 

\s

{\bf b.}  We prove that, for each $g \geqslant 4$ even, $C_{g,1} \cong C_{g,4}$ is the only Riemann surface of genus $g$ with a group of automorphisms of order $3g+3.$

\s

{\bf c.} We determine the automorphism groups of $C_{g,1}$, $C_{g,2}$ and $C_{g,4}$ and describe their corresponding actions. We deduce that they are non-hyperelliptic. 

\s

{\bf d.} We derive the fact that if $g$ is odd and not divisible by three, then there is no compact Riemann surface of genus $g$ possessing exactly $3g$ automorphisms.

\s

{\bf e.} We provide an isogeny decomposition of the Jacobian variety of $C_{g,1}$, $C_{g,2}$ and $C_{g,4}.$ In particular, we obtain that for each $g$ odd and not divisible by three, the Jacobian variety of $C_{g,2}$ and $C_{g,4}$ contain an elliptic curve.

\s

{\bf f.}  We show that if $g$ is even then there are infinitely many Riemann surfaces of genus $g$ with a group of automorphisms of order $3g$ and that the group is not necessarily cyclic. Furthermore, we give a  complete classification for the case $g=2p$ where $p$ prime.

\s

{\bf g.}  We show that if $g$ is odd then there are at least one Riemann surface of genus $g$ with a non-cyclic group of automorphisms of order $3g+3$. Furthermore, we give a  complete classification for the case $g=4p-1$ where $p$ prime.

\s

It is worth remarking that our results {\bf a} and {\bf b} above can be stated in terms of orientably-regular hypermaps (or dessin d'enfants); we shall briefly discuss that in \S\ref{hyper}.

\section{Preliminaries}
Let $X$ be a compact Riemann surface of genus $g \geqslant 2.$  By the classical uniformisation theorem, it is known that $X$ is isomorphic to a quotient of the form $\mathbb{H}/{\Gamma}$ where $\Gamma$ is a Fuchsian group, uniquely determined up to conjugation.

A finite group $G$ acts on $X$ if there is a group monomorphism $$\psi: G \to \mbox{Aut}(X),$$where $\mbox{Aut}(X)$ stands for the automorphism group of $X.$ In terms of Fuchsian groups, the group $G$ acts on $X \cong \mathbb{H}/{\Gamma}$ if and only if there is a Fuchsian group $\Delta$ and a group  epimorphism (ske, for short) \begin{equation*}\label{episs}\theta: \Delta \to G \mbox{ such that }  \mbox{ker}(\theta)=\Gamma.\end{equation*}

If the genus of $X/G \cong \mathbb{H}/\Delta$ is $\gamma$ and $k_1, \ldots, k_s$ are the branch indices in the canonical projection $\mathbb{H} \to \mathbb{H}/{\Delta}$, then the {\it signature} of $\Delta$ (and of the action of $G$) is the tuple \begin{equation} \label{sig} \sigma(\Delta)=(\gamma; k_1, \ldots, k_s).\end{equation}
In addition, in such a case the group $\Delta$ has the following canonical presentation  \begin{equation*}\label{prese}\langle \alpha_1, \ldots, \alpha_{\gamma}, \beta_1, \ldots, \beta_{\gamma}, x_1, \ldots , x_s: x_1^{k_1}=\cdots =x_s^{k_s}=\Pi_{i=1}^{\gamma}\alpha_i\beta_i\alpha_i^{-1}\beta_i^{-1} \Pi_{i=1}^s x_i=1\rangle.\end{equation*}We say that the action of $G$ is represented by $\theta$. With a slight abuse of notation, we also identify $\theta$ with the tuple of the images of the  
canonical generators of $\Delta.$

\s

The genus $g$ of $X$ and the signature of the action of $G$ are related by the equality$$2g-2=|G|(2\gamma-2+\Sigma_{j=1}^s(1-1/k_j)),$$known as the Riemann-Hurwitz formula.
\s

Let $\tau \in \psi(G) \cong G$ be a non-trivial automorphism of $X$. According to Macbeath \cite{McB73}, the number of fixed points of $\tau$ on $X$ is given by $$|\mbox{fix}(\tau)|= \mbox{N}_G(\langle \tau \rangle)\Sigma_{j=1}^s \epsilon_j(\tau)/k_j$$where $\epsilon_j(\tau)=1$ if and only if $\tau$ is conjugate to a power of $\theta(x_j).$ See also \cite{R07}.

\s

Two actions $\psi_i: G \to \mbox{Aut}(X)$  are {\it topologically equivalent} if there exist a orientation-preserving self-homeomorphism $h$ of $X$ and  $\omega \in \mbox{Aut}(G)$ such that
\begin{equation}\label{equivalentactions}
\psi_2(g) = h \psi_1(\omega(g)) h^{-1} \, \mbox{ for all } g \in G.
\end{equation}Each homeomorphism $h$ satisfying (\ref{equivalentactions}) induces an outer automorphism $h^*$ of $\Delta$. Let $\mathscr{B}$ denote the subgroup of $\mbox{Aut}(\Delta)$ consisting of those $h^*$. It is known that   $\theta_1, \theta_2 : \Delta \to G$ define actions that are topologically equivalent  if and only if $$\exists \omega \in \mbox{Aut}(G) \mbox{ and } \exists h^* \in \mathscr{B} \mbox{ such that }\theta_2 = \omega\circ\theta_1 \circ h^*.$$

We refer to \cite{B90}, \cite{H71} and \cite{McB} for more details.

\s

A {\it family} of compact Riemann surfaces of genus $g$ is the locus $$\mathcal{F}_g=\mathscr{M}_g(G, \sigma) \subset \mathscr{M}_g$$of points representing compact Riemann surfaces $S$ of genus $g$  with a group of automorphisms isomorphic to a given group $G$ acting with a given signature $\sigma$. It is known that each family consists of finitely many irreducible components that are in bijective correspondence with the topological classes of actions of $G$ with signature $\sigma.$ See \cite{B90}.

\s

If the signature of the action of $G$ on $X$ is \eqref{sig} then the complex dimension of the family is $3\gamma-3+s.$ Observe that zero-dimensional families correspond to Riemann surfaces with many automorphisms, namely, those surface that cannot be deformed in the moduli space keeping their symmetries. 

\s

We denote by $J_X$ the Jacobian variety of $X$. We recall that $J_X$ is a principally polarised abelian variety of dimension $g$ and that the Riemann surface $X$ is determined by $J_X.$ More precisely, the classical Torelli  theorem states that $$X_1 \cong X_2 \, \iff \, J_{X_1} \simeq J_{X_2},$$where $\simeq$ stands for isomorphism as principally polarised abelian varieties. It is known that each group of automorphisms $G$ of $X$ induces a group of automorphisms of $J_X$ that allows to provide an isogeny decomposition of $$J_X \sim B_1^{n_1} \times \cdots \times B_s^{n_s}$$ in terms of abelian subvarieties $B_i^{n_i} \subset J_X$ that are pairwise non $G$-isogenous. See \cite{CR06} and \cite{LR04} for more details. The dimension of each factor $B_i$ can be determined from the geometry of the action of $G$ on $X$. See \cite{R07}.

\s

Along the article, we denote the cyclic group of order $n$ by $\mathbb{Z}_n$ and the dihedral group of order $2n$ by $\mathbf{D}_{n}$. 

\section{Statement of the results}

\begin{theorem}\label{unicidad} Let $X$ be a compact Riemann surface of genus $g \geqslant 3.$ 

\s

Assume $g$ odd and different from 21. Then $X$ is endowed with a group of automorphisms of order $3g$ if and only if 
$$X \cong C_{g,2} \mbox{ if }g \not\equiv 2 \mbox{ mod }3 \, \mbox{ or } \, X \cong C_{g,4}   \mbox{ if } g \equiv 2 \mbox{ mod }3.$$

Assume $g$ even. Then $X$ is endowed with a group of automorphisms of order $3g+3$ if and only if  $$X \cong C_{g,1} \cong C_{g,3} \mbox{ and } g \not\equiv 2 \mbox{ mod } 3.$$
\end{theorem}

The existence of the Riemann surface in the theorem above also holds for $g=21$. However, the uniqueness is no longer true; we shall work out this sporadic case in \S\ref{rem:unicidad}.

\begin{theorem}\label{aut} Let $g \geqslant 3$ be an integer. If $g$ is odd then 
$$\mbox{Aut}(C_{g,2}) \cong \mathbb{Z}_{3g}, \,\, \mbox{Aut}(C_{g,2}) \cong \mathbb{Z}_3 \times \mathbf{D}_g \,\, \mbox{ and } \,\,  \mbox{Aut}(C_{g,4}) \cong \mathbb{Z}_3 \times \mathbf{D}_g$$for $g$ congruent to 0, 1 and 2 modulo 3, respectively. In the latter two cases, the automorphism group acts with signature $(0; 2,6,3g).$ If $g$ is even and $g \not\equiv 2 \mbox{ mod } 3$ then $$\mbox{Aut}(C_{g,1}) \cong \mathbb{Z}_{3g+3}.$$
\end{theorem}

\begin{corollary}\label{cor:aut} \mbox{}
\begin{enumerate}
\item The compact Riemann surfaces represented by the curves $C_{g,1}$, $C_{g,2}$ and $C_{g,4}$ are non-hyperelliptic, for each $g \geqslant 3$ as in the theorem above.
\s

\item If $g \geqslant 5$ is odd and $g \not\equiv 0 \mbox{ mod } 3$ then there does not exist a compact Riemann surface  genus $g$ with exactly $3g$ automorphisms. 
\end{enumerate}
\end{corollary}

\begin{theorem}\label{jaco} Let $g \geqslant 3$ be an integer and let $\varphi$ denote the Euler function.
\s

Assume $g$ odd and $g \not\equiv 0 \mbox{ mod }3$. Set $X=C_{g,2}$ if $g \equiv 1 \mbox{ mod } 3$ and $X=C_{g,4}$ if $g \equiv 2 \mbox{ mod } 3.$ The Jacobian  $J_{X}$ of $X$ decomposes, up to isogeny, as the product $$J_{X} \sim E \times \Pi_{d \in \mathcal{P}_g} A_{d}^2$$where $E$ and $A_d$ are abelian subvarieties of $J_{X}$ and $$\mathcal{P}_g=\{d \in \{1, \ldots, \tfrac{g-1}{2}\} :d \mbox{ divides }g\}.$$ The dimension of $A_d$ is $\tfrac{1}{2}\varphi(\tfrac{g}{d})$ and $E$ is an elliptic curve. In addition $$E \sim J_{X/\langle \varsigma \rangle}  \mbox{ and } \,\Pi_{d \in \mathcal{P}_g} A_{d} \sim  J_{X/\langle \iota \rangle}$$where $\varsigma$ and $\iota$ are automorphisms of $X$ of order $g$ and $2$ respectively.

\s

Assume $g$ odd and $g \equiv 0 \mbox{ mod }3.$ Set $X=C_{g,2}.$  The Jacobian $J_{X}$ of $X$  decomposes, up to isogeny, as the product$$J_{X} \sim  \Pi_{d \in \mathcal{Q}_g} B_{d}$$where $B_d$ is an abelian subvariety of $J_{X}$ of dimension $\tfrac{1}{2}\varphi(\tfrac{3g}{d})$ and $$\mathcal{Q}_g=\{d \in \{1, \ldots, 3g-1 \} : d \mbox{ divides } g \mbox{ and } d \not\equiv 0 \mbox{ mod }3\}.$$

Assume $g$ even and $g \not\equiv 2 \mbox{ mod } 3.$ Set $X=C_{g,1}$.
The Jacobian $J_{X}$ of $X$  decomposes, up to isogeny, as the product$$J_{X} \sim  \Pi_{d \in \mathcal{R}_g} B_{d}$$where $B_d$ is an abelian subvariety of $J_{X}$ of dimension $\varphi(\tfrac{g+1}{d})$ and $$\mathcal{R}_g=\{d \in \{1, \ldots, g \} : d \mbox{ divides } g+1 \}.$$
\end{theorem}

\begin{prop}\label{ct33g}
The elliptic curve $E$ in the theorem above is represented by $$w^3 = z^2-z.$$
\end{prop}

The following  results provide a complete classification of compact Riemann surfaces of genus $g$ with a group of automorphisms of order $3g$   and of order $3g+3,$ in a special case of $g$ even and $g$ odd, respectively.
\s

\begin{theorem} \label{gpar} Let $p \geqslant 7$ be a prime number, and let $X$ be a compact Riemann surface of genus $g=2p$.  If $X$ is endowed with a group of automorphisms of order $3g$ then one of the following statements holds.
\begin{enumerate}

\s

\item $X$ belongs to the complex one-dimensional family $\mathcal{F}_{g,1}$ of Riemann surfaces with a group of automorphisms $G$ isomorphic to $$\mathbf{D}_{3p} \mbox{ acting with signature }(0;2,2,3,3p).$$The members of $\mathcal{F}_{g,1}$ are trigonal and algebraically represented by$$y^{3}=x(x^p-1)(x^p-\eta)^{2} \mbox{ where } \eta \in \mathbb{C}- \{0,1\}.$$

\item $X$ belongs to the complex one-dimensional family $\mathcal{F}_{g,2}$ of Riemann surfaces with a group of automorphisms $G$ isomorphic to $$\mathbf{D}_{p} \times \mathbb{Z}_3 \mbox{ acting with signature }(0;2,2,3,3p).$$  The members of $\mathcal{F}_{g,2}$ are trigonal and algebraically represented by $$y^3=x^{\epsilon}(x^p-1)(x^p-\eta) \mbox{ where } \eta \in \mathbb{C}-\{0,1\},$$ and where $\epsilon=1$ if $p\equiv 2 \mbox{ mod }3$ and $\epsilon=2$ if $p\equiv 1 \mbox{ mod }3.$
\end{enumerate}

Both families are irreducible subvarieties of $\mathscr{M}_{2p}$ and, up to finitely many exceptions, the automorphism group of  their members is isomorphic to $G.$ 

\s

In addition, the unique compact Riemann surface of genus $g$ with an automorphism of order $3g$ belongs to $\mathcal{F}_{g,2}$, its  automorphism group  is isomorphic to $$\mathbf{D}_p \times \mathbb{Z}_6 \mbox{ and acts with signature } (0; 2, 6, 3g).$$
\end{theorem}

\begin{theorem} \label{gimpar} Let $p \geqslant 19$ be a prime number such that $p \equiv -1 \mbox{ mod } 4$, and let $X$ be a compact Riemann surface of genus $g=4p-1$.  If $X$ is endowed with a group of automorphisms of order $3g+3$ then one of the following statements holds.
\begin{enumerate}

\s
\item $X$ is isomorphic to $C_{g,1}$ if $g \equiv 0 \mbox{ mod }3$ or $C_{g,3}$  if  $g \equiv 1 \mbox{ mod }3$, and $G$ is cyclic.

\s

\item $X$ is isomorphic to the Riemann surface of genus $g$ with automorphism group isomorphic to $$\mathbb{Z}_p\rtimes_2\mathbb{Z}_{12}=\langle x, y: x^p=y^{12}=1, yxy^{-1}=x^{-1}\rangle$$ acting with signature $(0;4,12,3p).$ This Riemann surface  is $3p$-gonal  and algebraically represented by the curve $$y^{3p}=(x^2-1)(x^2+1)^{s}$$where $s \in \{2, \ldots, 3p-1\}$ is the unique integer satisfying $$
  s \equiv -1 \mbox{ mod } p  \mbox{ and }
  s \equiv 1  \mbox{ mod } 3.$$
\end{enumerate}
\end{theorem}

Theorem \ref{unicidad} and \ref{gpar} imply the following result.
\begin{theorem}\label{final}Let $g \geqslant 3$ be an integer. The compact Riemann surface of genus $g$ represented by the curve $$C_{g,2} \mbox{ if }g \not\equiv 2 \mbox{ mod }3 \, \mbox{ and } \, C_{g,4}   \mbox{ if } g \equiv 2 \mbox{ mod }3$$is uniquely determined by the the property of having a group of automorphisms of order $3g$ if and only if $g$ is odd and different from $21.$ 
\end{theorem}

\section{Orientably-regular hypermaps}\label{hyper}

By the classical Belyi theorem \cite{Belyi}, a compact Riemann surface $X$ of genus $g \geqslant 2$ can be defined over a number field if and only if $$\exists \beta :X \to \mathbb{P}^1 \mbox{ holomorphic with three critical values.}$$ In such a case, $X$ is uniformised by a Fuchsian group $\Gamma$ which is a finite index subgroup of a triangle Fuchsian group $$\Delta=\langle x,y : x^a=y^b=(xy)^c=1\rangle \mbox{ where } \tfrac{1}{a}+ \tfrac{1}{b}+ \tfrac{1}{c} <1.$$See, for instance, \cite{GG} and \cite{abc}.

\s

The map $\beta$ induces a so-called {\it dessin d’enfant} $$\mathscr{H}=\beta^{-1}([0,1])$$in $X$, namely, a finite, connected and bicoloured graph such that $X-\mathscr{H}$ is the union of finitely many topological discs. Equivalently, $\beta$ induces a (orientable) {\it hypermap}, with the preimages of $0$ providing the hypervertices, the preimages of $1$ the hyperedges, and the preimages of $\infty$ the hyperfaces. The {\it genus} of the hypermap is the genus of the supporting surface, and the hypermap is called {\it orientably-regular} if $\Gamma$ is a normal subgroup of $\Delta.$ In such a case, $\beta$ is a regular covering map associated to the action of $G\cong \Delta/\Gamma$ on $X$, and the automorphism group of the hypermap is isomorphic to $G$. 

\s

The problem of classifying hypermaps (and, in particular, {\it maps}) on (orientable and non-orientable) surfaces has been central in topological graph theory, and has been considered by several authors from different points of view. One of these points of view is by means of classifying according to its genus (or Euler characteristic) and to the order of its automorphism group. We refer, for instance, to the articles  \cite{CPS10}, \cite{CST10}, \cite{Aze05}, \cite{IJRC21} and \cite{Ma21}, and the references therein.

\s

In this terminology, our theorems \ref{unicidad} and \ref{final} can be restated  as follows.

\begin{theorem}
Let $g \geqslant 3$ be an integer.

\s

There exists a unique orientably-regular hypermap of genus $g$ with (orientation-preserving) automorphism group of order $3g$ if and only if $g$ is odd and different from $21.$ The hypermap is of type $(3,3g,3g)$  and is supported by the Riemann surface   $$C_{g,2} \mbox{ if }g \not\equiv 2 \mbox{ mod }3 \, \mbox{ and } \, C_{g,4}   \mbox{ if } g \equiv 2 \mbox{ mod }3.$$

\s
If $g$ is even and $g \not\equiv 2 \mbox{ mod } 3$ then there exists a unique orientably-regular hypermap of genus $g$ with (orientation-preservering)  automorphism group of order $3g+3.$ The hypermap is of type $(3,g+1,3g+3)$ and is supported by the  Riemann surface  $$C_{g,1} \cong C_{g,3}.$$
\end{theorem}
We refer to \cite{BCCI19} for the analogous situation for orientably-regular hypermaps of genus $g$ a group of automorphisms of order $4g$. We also refer to the database \cite{Conder} for several lists of orientably-regular maps and hypermaps of genus at most 101.

\section{Proof of Theorem \ref{unicidad}}

\subsection*{Riemann surfaces with $3g$ automorphisms}
Let $g \geqslant 3$ be an odd integer and let $G$ be a group of order $3g$.
 \subsection*{Case 1.1} Assume that $g \not\equiv 0 \mbox{ mod }3$. We recall that, by the famous odd-order theorem, the group $G$ is solvable. Since $g$ is not divisible by three,  Hall's theorem allows us to ensure that $G$ has a Hall subgroup $N$ such that $$[G:N]=3 \, \mbox{ and therefore } N \unlhd G.$$

If $H$ is a Sylow $3$-subgroup of $G$ then, by the Schur-Zassenhaus theorem,  one has  $$G \cong N \rtimes H \, \mbox{ where } \, H \cong \mathbb{Z}_3.$$

Let $X$ be a compact Riemann surface of genus $g$ endowed with a group of automorphisms isomorphic to $G$, and consider the quotient Riemann surfaces  $$Y:=X/N \, \mbox{ and } \, Z:=X/H$$given by the action of $N$ and $H$ on $X$.  We denote their genera by $g_Y$ and $g_Z$ respectively. An easy consequence of the Riemann-Hurwitz formula is that $g_Y$ equals either $0$ or $1$.

\s

{\it Claim.} If $\tau \in G$ and $\langle \tau \rangle=H$ then $\tau$ has fixed points. 

\s

The fact that each element of $G$ of order three is conjugate to $\tau$ or to $\tau^2$ implies that if $\tau$ does not have fixed points then there is no point in $X$ with $G$-stabiliser of order divisible by three. It then follows that the automorphism $\bar{\tau}$ of $Y$ induced by $\tau$ does not have fixed points in $Y$. Hence the signature of the action of $N$ on $X$ is  $$(1; e_1, e_1, e_1, \ldots, e_l, e_l, e_l) \, \mbox{ for some $e_i \geqslant 5$ and }\, l \geqslant 1.$$However, the Riemann-Hurwitz formula applied to the covering $X \to Y$ implies that
$$
2g-2 \geqslant g[3(1-\tfrac{1}{e_1})] \geqslant \tfrac{12}{5}g \implies g \leqslant -5
$$and the claim follows.

\s

Observe that the number $v$ of points of $Y$ fixed by $\bar{\tau}$ is $2$ or $3$ according to whether $g_Y$ equals $0$ and $1$, and that in both cases $X/G \cong \mathbb{P}^1.$ Assume that there is no element in $N$ which fixes a point that is also fixed by $\tau$. Then, the signature of the action of $G$ on $X$ has the form $$(0; 3,3,e_1, \ldots, e_l) \, \mbox{ or } \, (0; 3,3,3,e_1, \ldots, e_l) \mbox{ where }l \geqslant 1,$$for some $e_j \geqslant 4,$ according to $g_Y$ equals $0$ or $1.$ A simple computation shows that both cases contradict the Riemann-Hurwitz formula. It then follows that there is a fixed point of $\tau$ which is also fixed by some element of $N.$ 

\s

Note that if $g_Y=0$ then exactly one or two $H$-short orbits are also $N$-short orbits. In other words, the signature of the action of $G$ on $X$ has the form $$(0; 3,3e_1, e_2, \ldots, e_l) \, \mbox{ or } \, (0; 3e_1, 3e_2, e_3, \ldots, e_l)$$and again this contradicts the Riemann-Hurwitz formula. We obtain that $g_{Y}=1$ and the signature of the action of $G$ has the form
 $$(0; 3e_1, 3e_2, 3e_3, e_4, \ldots, e_l),(0; 3, 3e_2, 3e_3, e_4, \ldots, e_l) \mbox{ or }(0; 3, 3, 3e_3, e_4, \ldots, e_l),$$where $5 \leqslant e_i \leqslant g$ divides $g.$ The former and latter cases are ruled out similarly as before. A computation shows that the signature of the action is necessarily of the form $$(0; 3,3e_1, 3e_2) \mbox{ where }\tfrac{1}{e_1}+\tfrac{1}{e_2}=\tfrac{2}{g}.$$
 
 Since the unique solution of the previous equation is $e_2=e_3=g$, we conclude that the signature of the action of $G$ on $X$ is $(0; 3,3g,3g).$ Consequently, the group $G$ is cyclic and the uniqueness results in \cite{N84} ensure that $$X \cong C_{g,2} \mbox{ if } g \equiv 1 \mbox{ mod }3 \mbox{ or }X \cong C_{g,4} \mbox{ if }g \equiv 2 \mbox{ mod }3.$$

 \subsection*{Case 1.2} Assume that $g \equiv 0 \mbox{ mod }3$. Since it is well-known that there is only one Riemann surface of genus three with a group of automorphisms of order nine (see, for instance, \cite{B91}, \cite{Conder} and \cite{Kn}), we can assume $g \geqslant 9.$ 
 
 \s
 
 Let $n \geqslant 3$ be an odd integer, and let  $X$ be a compact Riemann surface of genus $g=3n$ endowed with a group of automorphisms isomorphic to $G$ (of order $9n$). Assume that the signature of the action of $G$ on $X$ is $(\gamma; m_1, \ldots, m_r)$ and let $$B=\Sigma_{i=1}^r(1-1/m_i).$$Observe that if $\gamma \geqslant 2$ then the Riemann-Hurwitz formula implies that $$6n-2 \geqslant 9n(2+B) =18n+9nB,$$ which is impossible. Besides, if $\gamma=1$ then the fact that $m_i \geqslant 3$ shows that $6n-2 \geqslant 9nB \geqslant 6nr$ and therefore $r=0,$ a contradiction. It follows that $\gamma=0$ and consequently $$r \leqslant 4-\tfrac{1}{3n} \mbox{ showing that } r=3.$$ Thus, the signature of the action of $G$ is $(0; m_1, m_2, m_3)$ where \begin{equation}\label{tel}\tfrac{1}{m_1}+\tfrac{1}{m_2}+\tfrac{1}{m_3}=\tfrac{1}{3}+\tfrac{2}{9n}.\end{equation}
 
 \s
 
{\it Claim.} The solutions of  \eqref{tel} are $m_1=3, m_2=m_3=9n$ for each $n \geqslant 3$, and six sporadic cases given in the following table
 \s
 \begin{center}
\begin{tabular}{|c|c|c|c|c||c|c|c|c|c|}  
\hline
case & $m_1$ & $m_2$ &  $m_3$  &  $n$ & case  & $m_1$ & $m_2$ &  $m_3$  &  $n$\\ \hline 
{\bf a} & 5 & 9  & 15 & 5 & {\bf d} & 5 & 9  & 39 & 65  \\ 
{\bf b} & 5 & 9  & 27 & 15 & {\bf e} & 5 & 9  & 43 & 215  \\ 
{\bf c} & 5 & 9  & 35 & 35 &{\bf f} & 7 & 9  & 9 & 7  \\ 
\hline
\end{tabular}
\end{center}

\s

Without loss of generality, we can assume $m_1 \leqslant m_2 \leqslant m_3.$ Let $u$ be the number of periods that are equal to 3. Clearly $u \neq 3$. If $u=2$ then $m_3=9n/(2-3n) <0$, and if $u=1$ then $m_1=3$ and $m_2=m_3=9n.$ We assume that $u=0.$ Observe that, if $m_1, m_2, m_3$ were greater than or equal to 9, then $n <0$. Thus, $m_1=5$ or $m_1=7$.

\s

Assume $m_1=5.$ If $m_2, m_3 \geqslant 15$ then $n <0,$ and then $m_2 \in \{5, 7,9,11,13\}.$

\begin{enumerate}
\item If $m_2=9$ then $m_3=\tfrac{45n}{10+n}$ and $n$ is divisible by 5. Hence $$(n,m_3) \in \{(5,15), (15,27), (35,35), (65,39), (215,43)\}$$ 
\item The remaining possibilities for $m_2$ do not give rise to any solution of  \eqref{tel}. For instance, if $m_2=7$ then $m_3=\tfrac{315n}{70-3n}$  and $n$ must be divisible by 35; this is not possible. The other cases are analogous.
\end{enumerate}

Assume $m_1=7.$ If $m_2, m_3 \geqslant 11$ then $n <0,$ and then $m_2 \in \{ 7,9\}.$ Analogously as before, one sees that the only solution of  \eqref{tel} is $m_2=m_3=9$ and $n=7.$ This proves the claim.

\s

The cases {\bf a} and {\bf b} are not realised; see for instance \cite{Conder}. The cases {\bf c}, {\bf d} and {\bf e} are not realised neither. In fact, assume that there exists a compact Riemann surface $Y$ of genus $105$ endowed with a group of automorphisms $G$ of order $315.$ According to the table above, we also assume that $G$ acts on $Y$ with signature $(0;5,9,35).$ Clearly $G$ cannot be abelian since the product of an element of order 5 and one of order 9 must have order 45. Now, since $G$ has an element of order 9 and $|G| = 3^2\times 5\times 7$, one sees that each Sylow $3$-subgroup of $G$ is a cyclic group of order $9$. Among the non-cyclic groups of order 315, there is exactly one with cyclic Sylow $3$-subgroups, as it can be verified by using Magma (the group is labelled as $(315,1)$ in the SmallGroup library). This group has a normal Sylow $5$-subgroup $K$ of order $5$. The normality of $K$ ensures that it acts on $Y$ fixing 72 points. In other words,  the regular covering map $Y \to Y/K$ given by the action of $K$ on $Y$ ramifies over 72 points and hence $$208=10(\gamma-1)+288, \mbox{ where $\gamma$ is the genus of $Y/K$.}$$ This is clearly impossible. The cases {\bf d} and {\bf e} are ruled out analogously.

We then conclude that if $g \neq 21$ then $G$ is cyclic and, by the uniqueness results in \cite{N84}, $X$ is isomorphic to $C_{g,2}$.

\subsection*{Riemann surfaces with $3g+3$ automorphisms} Let $g \geqslant 4$ be even, write $n:=g+1$ and let $G$ be a group of order $3n=3g+3$. 
 \subsection*{Case 2.1} Assume that $n \not\equiv 0 \mbox{ mod }3$. Since $G$ is solvable and $n$ is not divisible by three,  there is a normal subgroup $N$ of $G$ of order $n$ in such a way that $$G \cong N \rtimes H \, \mbox{ where } \, H \cong \mathbb{Z}_3.$$

Let $X$ be a compact Riemann surface of genus $n-1$ endowed with a group of automorphisms isomorphic to $G$. Consider the quotient Riemann surface  $Y:=X/N,$ whose genus $g_Y$ equals either $0$ or $1$. By arguing as in Case 1.1, one sees that if $\langle \tau \rangle=H$ then $\tau$ has fixed points, and that among such fixed points there is at least one which is also fixed by some element of $N.$ This fact, together with the assumption $g_Y=1$ imply that the signature of the action of $G$ on $X$ is$$(0; 3e_1, 3e_2, 3e_3, e_4, \ldots, e_l),(0; 3, 3e_2, 3e_3, e_4, \ldots, e_l) \mbox{ or }(0; 3, 3, 3e_3, e_4, \ldots, e_l),$$where $e_i \geqslant 5$ and $l \geqslant 4.$ However, these signatures contradict the Riemann-Hurwitz formula. It follows that $g_{Y}=0$ and hence the signature of the action of $G$ is$$(0; 3,3e_2, e_3, \ldots, e_l) \, \mbox{ or } \, (0; 3e_1, 3e_2, e_3, \ldots, e_l)$$where $5 \leqslant e_i \leqslant n$ and $l \geqslant 3.$ The latter case can be ruled out similarly as before. A computation shows that the signature of the action is necessarily of the form $$(0; 3,3e_2, e_3) \mbox{ where }\tfrac{1}{e_2}+\tfrac{3}{e_3}=\tfrac{4}{n}.$$Since the unique solution of the previous equation is $e_2=e_3=n$ we conclude that $G$ is cyclic. Now, the uniqueness results of \cite{N84} implies that that $$X \cong C_{g,1} \mbox{ if } g \equiv 0 \mbox{ mod }3 \mbox{ and }X \cong C_{g,3} \mbox{ if }g \equiv 1 \mbox{ mod }3.$$Finally, it is easily seen that the two curves above are both isomorphic to the curve $$y^3=x^{g+1}-1,$$and the proof is complete.

 \subsection*{Case 2.2} Assume that $n \equiv 0 \mbox{ mod }3$ and write $n=3m$ for some $m \geqslant 3$ odd. Let $X$ be a  Riemann surface of genus $g=3m-1$ endowed with a group of automorphisms $G$ of order $9m$.  A straightforward computation shows that necessarily the genus of $X/G$ is zero and that the corresponding regular covering map ramifies over exactly three values. Therefore, the signature of the action of $G$ is $(0; m_1, m_2, m_3)$ where $$\tfrac{1}{m_1}+\tfrac{1}{m_2}+\tfrac{1}{m_3}=\tfrac{1}{3}+\tfrac{4}{9m}$$

If $v$ is the number of periods that are equal to 3 then clearly $v \neq 2, 3$. We observe that $v \neq 1.$ Indeed, in such a case we may assume $m_1=3$ and therefore $\tfrac{1}{m_2}+\tfrac{1}{m_3}=\tfrac{4}{9m}.$ The fact that $m$ is odd implies that $m_2=3m$ and $m_3=9m,$  and therefore $G$ is cyclic. However, in a cyclic group the product of an element of order 3 and one of order $3m$ cannot be $9m$. We may then assume that $v=0.$ Observe that if $m_1, m_2, m_3 \geqslant 9$ then $m<0$, and therefore there are two cases to consider: $m_1=5$ or $m_1=7.$

\s

Assume that $m_1=5.$ If $m_2, m_3 \geqslant 15$ then $m <0,$ and therefore  $m_2 \in \{5, 7,9,11,13\}.$

\begin{enumerate}
\item If $m_2=5$ or $m_2=7$ then $m_3$ equals to $\tfrac{45m}{20-3m}$ or $\tfrac{315m}{140-3m}$ respectively. These cases only yield to the following signatures:  $$(0; 5,5,45) \mbox{ with } |G|=45 \, \mbox{ and } (0; 5,7,315) \mbox{ with } |G|=315.$$ In both cases $G$ is cyclic, but such signatures are not realised by  cyclic groups.

\s

\item If $m_2=9$, then $m_3=\tfrac{45m}{m+20}$ and therefore $$(m,m_3) \in \{(5,9), (25,25), (55,33), (205,41)\}.$$ 
The former two cases are not realised, as shown in \cite{Conder}. The latter two cases are not realised neither. In fact, assume that $X$ is a Riemann surface of genus $164$ (resp. $614$) equipped with a group of automorphisms $G$ of order $495$ (resp. $1845$) acting with signature $(0;5,9, 33)$ (resp. $(0;5,9,41)$). Note that there are exactly four groups $G$ of order $495$ (resp. $1845$). The signature of the action allows us to exclude the abelian ones. Now, since a Sylow $3$-subgroup $H$ of $G$ must be cyclic of order $9$, the group $G$ must be isomorphic to SmallGroup(495,1) (resp. SmallGroup(1845,1)). In addition, as $H$ is normal in $G$, we have that $H$ must fix $55$ (resp. $205$) points. However, an easy application of the Riemann-Hurwitz formula to the regular covering map $X\to X/H$ gives a contradiction.

\s 
\item If $m_2=11$ or $m_2=13$ then $m_3$ equals to $\tfrac{495m}{220+21m}$ or $\tfrac{585m}{260+33m}$ respectively. Both cases do not give rise to any solution. \end{enumerate}

Assume $m_1=7.$ If $m_2, m_3 \geqslant 11$ then $n <0.$ Thus we can assume that $m_2 \in \{7,9\}.$ By proceeding analogously as before, one sees that the only two possibilities are  $$(0; 7,7,9), m=7, |G|=63 \mbox{ and } (0; 7,9,9), m=14, |G|=126.$$Both cases are not realised, as shown in \cite{Conder}.

\subsection*{Remarks on Theorem \ref{unicidad}} \label{rem:unicidad} \mbox{}

\s

{\it The exceptional case in Theorem \ref{unicidad}.} The group $$\mathbb{Z}_7 \rtimes_3 \mathbb{Z}_9=\langle a,b : a^7=b^9=1, bab^{-1}=a^4\rangle$$is non-cyclic and has order $63$. The tuple $(a,b,(ab)^{-1})$ is a ske of $G$ of signature $s=(0; 7,9,9)$ and hence there is a  Riemann surface $Y$ of genus $21$ with an action of $G$ with signature $s$. Since the automorphism group of $Y$ is isomorphic to $G$ (see, for instance,  \cite{Conder}), Theorem \ref{aut} implies that $Y$ is not isomorphic to $C_{21,2}.$ This is the unique exception in Theorem \ref{unicidad}.  

\s

{\it The case $g$ prime in Case 1.1.}  If $g$ is prime, then it was shown in \cite[Proposition 11]{AS} that a cyclic group of automorphisms $N$ of order $g$  of a compact Riemann surface $X$ of genus $g$ fixes two points and that  $X/N$ is elliptic. This fact, coupled with  the existence of an automorphism $\tau$ of order three normalising $H$ in Case 1.1 in the proof of Theorem \ref{unicidad}, immediately allow to conclude that the group generated by $\tau$ and a generator $\sigma$ of $N$ is cyclic (as $\tau$ must act on the points that are fixed by $N$) and hence, the classification follows.

\section{Proof of Theorem \ref{aut}}
\subsection*{Case 1}  Let $g \geqslant 5$ be an odd integer such that $g \not\equiv 0 \mbox{ mod }3$, and consider  $$G=\mathbb{Z}_g \times \mathbb{Z}_3=\langle a, b: a^g=b^3=[a,b]=1\rangle.$$

The elements of order three are $b$ and $b^{-1}$, whereas the ones of order $3g$  are $$a^kb^{\pm 1} \mbox{ where } 1 \leqslant k \leqslant g-1 \mbox{ such that }(k,g)=1.$$
 
We write $X=C_{g,2}$ if $g \equiv 1 \mbox{ mod }3$ and $X=C_{g,4}$ if $g \equiv 2 \mbox{ mod }3.$ Since the signature of the action of $G$ on $X$ is $(0; 3,3g,3g),$ one has that such an action is represented by any ske of the form $$(b^{\pm 1}, a^kb^{\pm 1}, a^lb^{\pm 1}) \mbox{ for some }k,l \mbox{ satisfying }(k,g)=(l,g)=1.$$Without loss of generality, we choose $\theta:=(b,ab,a^{-1}b)$ representing the action. 

\s

Consider the group $$\bar{G}=\langle t,r,s: t^3=r^g=s^2=(sr)^2=[t,r]=[t,s]=1\rangle \cong \mathbb{Z}_3 \times \mathbf{D}_g$$of order $6g.$ The tuple\begin{equation*}\label{skegrande}\bar{\theta}=(s,tsr,(tr)^{-1})\end{equation*}is a ske representing an action of $\bar{G}$ on a compact Riemann surface $\bar{X}$ with signature $(0; 2,6,3g).$ The genus of $\bar{X}$ is $g.$  As $\langle t,r \rangle$ is a subgroup of $\bar{G}$ isomorphic to $\mathbb{Z}_3 \times \mathbb{Z}_g$, the uniqueness of $X$ allows us to conclude that $X \cong \bar{X}.$ Finally, the fact that the signature $(0; 2,6,3g)$ is maximal  for each $g$ (see \cite{Sing72}) implies that $X$ does not admit more automorphisms, and hence the automorphism group of $X$ is isomorphic to $\mathbb{Z}_3 \times \mathbf{D}_g$.

\subsection*{Case 2} 
Let $g \geqslant 3$ be an odd integer such that $g \equiv 0 \mbox{ mod }3$. We write $g=3n$ where $n \geqslant 1$ is odd.  We write $$G=\mathbb{Z}_{9n}=\langle a: a^{9n}=1\rangle.$$The elements of order three of $G$ are $a^{\pm 3n}$ whereas the ones of order $9n$ are $$a^k \mbox{ where } 1 \leqslant k \leqslant 9n-1 \mbox{ such that }(k,9n)=1.$$

It follows that the action of $G$ on $X=C_{g,2}$ is represented by one of the ske $$(a^{\pm 3n}, a^k, a^{\pm 6n-k}).$$Indeed, the uniqueness of the Riemann surface allows us to choose $\theta=(a^{3n}, a, a^{6n-1})$ representing the action.
 
\s

Assume the automorphism group $G'$ of $X$ to have order strictly greater than $9n$. Then, following \cite{Sing72}, the group $G'$  has order $18n$ and acts on $X$ with signature $(0; 2,6,9n).$ Let $x,y,z$ be elements of $G'$ such that $\Theta=(x,y,z)$ is a ske representing the action of $G'$ on $X$. Since $G'$ has an index two cyclic subgroup and since $n$ is odd, the Schur-Zassenhaus theorem says that $$G' \cong \mathbb{Z}_{9n} \rtimes \mathbb{Z}_2=\langle A, B : A^{9n}=B^2=1, BAB=A^{s}\rangle$$for some $s \in \{1, \ldots, 9n-1\}$ such that $s^2\equiv 1 \mbox{ mod } 9n.$ Observe that $s \neq 1, 9n-1.$ Indeed, in the former case $G'$ is cyclic and therefore the fact that $x$ is an involution and $y$ has order six contradicts that fact that $z=(xy)^{-1}$ has order $9n.$ The latter case is impossible since the dihedral group does not have elements of order six.

Observe that the involutions of $G'$ are the elements of the form $$A^uB \mbox{ where } u(s+1) \equiv 0 \mbox{ mod }9n,$$and the elements of order six are those of the form $$A^vB \mbox{ where } 3v(s+1) \equiv 0 \mbox{ mod }9n \mbox{ and } v(s+1) \not\equiv 0 \mbox{ mod } 9n.$$
 
It then follows that $g:=z^{-1}=xy$, whose order is $9n,$ must have the form  $$g=A^uB \cdot A^vB=A^{u+sv}$$

As $(s,9n)=1$ one sees that $3s \not\equiv 0 \mbox{ mod }9n.$ Notice that $$g^{3s}=A^{3s(u+sv)}=A^{3us+3v}=A^{-3u-3vs}=g^{-3}$$and therefore $g^{3(s+1)}=1.$ It follows that $3(s+1) \equiv 0 \mbox{ mod }9n$ and therefore $s+1=3nl$ for some $l \in \mathbb{Z}.$ The fact that $s^2\equiv 1 \mbox{ mod }9n$ implies that $l=3l'$ for some $l' \in \mathbb{Z}$ and thereby $$s+1=9nl' \mbox{ showing that } s \equiv -1 \mbox{ mod } 9n,$$a contradiction. It then follows that the automorphism group of $X$ is isomorphic to $G$.

\subsection*{Case 3} Assume that $g$ is even and that $g \not\equiv 2 \mbox{ mod } 3.$ Since $C_{g,1}$ and $C_{g,3}$ are isomorphic to the curve $y^3=x^{g+1}-1$, the result  follows  from \cite[Theorem 1]{K98} where the automorphism group of curves of the form $x^n+y^m=1$ are studied.

\subsection*{Proof of Corollary \ref{cor:aut}}
The first statement follows from from the fact that $\mbox{Aut}(C_{g,2})$ and $\mbox{Aut}(C_{g,4})$ does not have central involutions (in the case of $C_{g,1}$  the result is obvious). The second statement follows directly from Theorems \ref{unicidad} and  \ref{aut}.

\section{Proof of Theorem \ref{jaco}}
We denote the Euler function by $\varphi$ and write $\omega_n=\mbox{exp}(2 \pi i /n)$ for each $n \geqslant 2$ integer.
\subsection*{Case 1}  Let $g \geqslant 5$ be an odd integer such that $g \not\equiv 0 \mbox{ mod }3$. We consider the following complex irreducible representations of $$\bar{G}=\langle t,r,s: t^3=r^g=s^2=(sr)^2=[t,r]=[t,s]=1\rangle \cong \mathbb{Z}_3 \times \mathbf{D}_g.$$

\begin{enumerate}
\item For each $k \in \{1, \ldots, \tfrac{g-1}{2}\}$, the representation of degree two given by $$V_k: r \mapsto \mbox{diag}(\omega_g^k, \bar{\omega}_g^k), s \mapsto \begin{psmallmatrix}0 & 1\\1 & 0\end{psmallmatrix}, t \mapsto \mbox{diag}(\omega_3, \omega_3).$$ 
\item The representation of degree one given by $V_0 :  r \mapsto 1, s \mapsto -1, t \mapsto \omega_3$.
\end{enumerate}

\s

Let $d \in \{1, \ldots, \tfrac{g-1}{2}\}$ such that $d$ divides $g.$ The character field of $V_{d}$ $$\mathbb{Q}(\omega_{g/d}+\bar{\omega}_{g/d}, \omega_3) \mbox{ has degree }\varphi(\tfrac{g}{d})\mbox{ over }\mathbb{Q}.$$

If we denote by $\mathscr{G}_d$ the Galois group corresponding to this field extension, then $$\mathcal{W}_d:=\oplus_{\sigma \in \mathscr{G}_d}V_{d}^{\sigma}$$is a rational irreducible representation of $\bar{G}$ (see, for instance, \cite[\S12.2]{Serre}). Indeed, if $\mathcal{W}_0$ is the direct sum of $V_0$ and its complex-conjugate representation, then the set$$\{\mathcal{W}_0, \mathcal{W}_d\}_{d \in \mathcal{P}_g} \mbox{ where } \mathcal{P}_g:=\{ d \in \{1, \ldots, \tfrac{g-1}{2}\}: d \mbox{ divides }g\}$$is formed by pairwise non-equivalent rational irreducible representation of $\bar{G}$.

\s

We write $X=C_{g,2}$ if $g \equiv 1 \mbox{ mod }3$ and $X=C_{g,4}$ if $g \equiv 2 \mbox{ mod }3.$ We recall that, as proved in Theorem \ref{aut},  the automorphism group of $X$ is isomorphic to $\bar{G}$ and acts on it with signature $(0; 2,6,3g).$ Thus, following the results of \cite{LR04}, there exists an abelian subvariety $T$ of $J_X$ such that $J_X$ is isogenous to the product $$A_0 \times \Pi_{d \in \mathcal{P}_g} A_{d}^2 \times T,$$where $A_j$ is an abelian subvariety of $J_X$ associated to the representation $\mathcal{W}_j$.

\s

In order to compute the dimension of the factors in the previous isogeny decomposition of $J_X$, we have to take into consideration the geometry of the action of $\bar{G}$ on $X$, represented by the ske \begin{equation*}\label{skeesc}\bar{\theta}=(s,tsr,(tr)^{-1}).\end{equation*} Observe that for each $d \in \mathcal{P}_g$ one has that $$\dim V_{d}^{\langle s \rangle}= 1 \mbox{ and } \dim V_{d}^{\langle tsr^{-1} \rangle}=\dim V_{d}^{\langle tr \rangle}=0$$(dimension of fixed subspaces) and consequently, by \cite[Theorem 5.12]{R07},$$\dim(A_{d})=\varphi(\tfrac{g}{d})[-2+\tfrac{1}{2}(6-1)]=\tfrac{1}{2}\varphi(\tfrac{g}{d}).$$  Similarly one can see that the dimension of $E:=A_0$ is 1. It follows  that $$\mbox{dim}(T)=g-1-2\Sigma_{d \in \mathcal{P}_g}\tfrac{1}{2}\varphi(\tfrac{g}{d})=0;$$ whence we conclude that $$J_X \sim E \times \Pi_{d \in \mathcal{P}_g} A^2_{d}.$$

\s

Let $\varsigma$ and $\iota$ be automorphisms of $X$ of order $g$ and $2$ respectively. Then, following \cite[Proposition 5.2]{CR06}, we have that $$J_{X/\langle \iota \rangle} \sim E^{k_0} \times \Pi_{d \in \mathcal{P}_g}A_{d}^{k_d} \mbox{ where } k_d=\dim V_{d}^{\langle \iota \rangle}.$$ The fact that $\iota$ is conjugate to $s$ shows that $k_0=0$ and $k_d=1$ for $d \neq 0.$ Thus, $$J_{X/\langle \iota \rangle} \cong J_{X/\langle s \rangle} \sim \Pi_{d \in \mathcal{P}_g}A_{d}.$$In a very similar way we obtain that $J_{X_{\langle \varsigma\rangle}} \cong J_{X_{\langle r \rangle}} \sim E.$ 

\begin{rema}
Our decomposition of $J_X$ could be refined if we had extra information on  $Y:=X/\langle s \rangle$, as its Jacobian  $J_Y$ is a factor  in the isogeny decomposition $$J_X \sim E \times J_{Y}^2.$$
The fact that  $\mathbf{t}$ and $\mathbf{s}$ commute implies that  $Y$ is trigonal. However, no further information can be gathered on $Y$ by Galois theory, so we cannot know \emph{a priori} if $Y$ has further automorphisms. A way to look into this would be to find an explicit algebraic model for $Y$. Observe that, as $X$ is a cyclic six-fold covering of $\mathbb{P}^1$, there exist $u,v \in \mathbb{C}(x)$ such that $u^6 = f(v)$ for some polynomial $f \in \mathbb{C}[x,y]$. Consequently, by taking $t = u^2$, we see that an equation for $Y$ is of the form $t^3 = f(v).$
\end{rema}

\subsection*{Case 2}  Let $g \geqslant 3$ be an odd integer such that $g \equiv 0 \mbox{ mod }3$. We recall that the complex irreducible representations of $G=\langle a: a^{3g}=1\rangle$ are$$\chi_{d}: a \mapsto \omega_{3g}^d \mbox{ where } d\in \{1 \ldots, 3g\},$$ and also that the action of $G$ on $X=C_{g,2}$ has signature $(0; 3,3g,3g)$ and is represented by the ske $\theta=(a^{g}, a, a^{2g-1})$. 

\s

We denote by $B_d$ the abelian subvariety of $J_X$  associated to $\chi_d.$ Since $\chi_{3g}$ is the trivial representation, one sees that $B_{3g}$ is zero-dimensional. In order to compute the dimension of the remaining $B_d$, observe that the dimension of the fixed subspace of $\chi_d$ by $\langle a \rangle$ equals zero for each $d$. Now, we apply \cite[Theorem 5.12]{R07} to see that  $$\dim (B_{d})=k_d[-1+\tfrac{1}{2}(3-\epsilon_d)] \, \mbox{ where } \epsilon_d=\dim \chi_d^{\langle a^g\rangle}$$and $k_d$ is the degree of the character field of $\chi_d$ over $\mathbb{Q}.$ Notice that $$\chi_d(a^g)=\omega_{3g}^{dg}=1 \, \iff \, dg \equiv 0 \mbox{ mod }3g$$and this, in turn, is equivalent to $\dim (B_{d})=0$. Hence, $\dim(B_d)\neq 0$ if and only if $$d \in \mathcal{Q}_g:=\{d \in \{1, \ldots, 3g -1\} : d \mbox{ divides } g \mbox{ and } d \not\equiv 0 \mbox{ mod }3\}.$$ Thus, following \cite{LR04}, we conclude that  $$J_X \sim \Pi_{d \in \mathcal{Q}_g} B_{d}.$$ Since the character field of $\chi_d$  is $\mathbb{Q}(\omega_{3g/d})$ for each $d \in \mathcal{Q}_g$, we obtain that $$ \dim B_{d}=\tfrac{1}{2}k_d=\tfrac{1}{2}\varphi(\tfrac{3g}{d}).$$

\begin{rema}
As a consequence of \cite[Proposition 5.2]{CR06}, one sees that if $m \in \mathcal{Q}_g$  then $$J_{X/H_m} \sim \Pi_{d \in \mathcal{Q}_g, m | d}B_{d} \mbox{ where } H_m=\langle a^{3g/m}\rangle \cong \mathbb{Z}_m.$$
\end{rema}

{\bf Case 3.} Assume $g$ even and $g \not\equiv 2 \mbox{ mod } 3$. Consider the pairwise non-equivalent complex irreducible representations of $$G=\langle a, b : a^3=b^{g+1}=[a,b]=1\rangle \cong \mathbb{Z}_3 \times \mathbb{Z}_{g+1}$$given by $$\rho_d : a \mapsto \omega_3, \, b \mapsto \omega_{g+1}^d$$for each $d \in \{1, \ldots, g\}$ such that $d$ divides $g+1.$ The Galois group associated to the character field of $\rho_d$ over $\mathbb{Q}$ has order $2\varphi(\tfrac{g+1}{d})$. Since the set $\{\rho_d\}$ consists of pairwise Galois non-equivalent representations, we have that $J_X$ of $X=C_{g,1}$ decomposes as $$J_{X} \sim  \Pi_{d \in \mathcal{R}_g} B_{d} \times T$$where $B_d$ is an abelian subvariety of $J_{X}$ associated to $\rho_d$, $T$ is an abelian subvariety of $J_X$  and $$\mathcal{R}_g=\{d \in \{1, \ldots, g \} : d \mbox{ divides } g+1 \}.$$

After taking into consideration that the action of $G$ on $X$ is represented by the ske $(a,b,(ab)^{-1})$, similarly as done before, we apply \cite[Theorem 5.12]{R07} to see that  $$\dim (B_{d})=2\varphi(\tfrac{g+1}{d})[-1+\tfrac{1}{2}(3-0)]=\varphi(\tfrac{g+1}{d}).$$Since $\Sigma_{d \in \mathcal{R}_g} \dim(B_{d})=g$ we conclude that $\dim(T)=0$ and the result follows.

\begin{rema}
As in the previous case, by \cite[Proposition 5.2]{CR06},  if $m \in \mathcal{R}_g$  then $$J_{X/H_m} \sim \Pi_{d \in \mathcal{R}_g, m | d}B_{d} \mbox{ where } H_m=\langle b^{(g+1)/m}\rangle \cong \mathbb{Z}_m.$$
\end{rema}

\subsection*{Proof of Proposition \ref{ct33g}}

Assume that $g$ is odd and $g \not\equiv 0 \mbox{ mod }3$. We write $$X=C_{g,2}: y^3=x(x^g-1) \mbox{ and } X=C_{g,4}: y^3=x^2(x^g-1)$$ if $g \equiv 1 \mbox{ mod }3$ and $g \equiv 2 \mbox { mod }3,$ respectively.  We recall that, in both cases, the automorphism group of $X$ is isomorphic to$$\langle t,r,s: t^3=r^g=s^2=(sr)^2=[t,r]=[t,s]=1\rangle \cong \mathbb{Z}_3 \times \mathbf{D}_g$$and that the quotient Riemann surface $E$ given by the action of $\langle r \rangle$ on $X$ has genus one. Consider the compact Riemann surface $Z$ represented by the algebraic curve$$C_{g,24}: y^3 = x^{2g}-x^g$$for each $g.$ Since $Z$ is a $3$-fold regular cover of $\mathbb{P}^1$ ramified over $g+2$ values, its genus is $g.$ We claim that $C_{g,2}$ and $C_{g,4}$ are isomorphic  to $C_{g,24}$. Indeed, the maps$$\mathbf{t}(x,y) = (x,\omega_3y) \mbox{ and } \mathbf{r}(x,y) = (\omega_gx,y)$$are automorphisms of $C_{g,24}$ of order 3 and $g$ that commute. It follows that $Z$ has a cyclic group of automorphisms of order $3g$ and the claim follows by the uniqueness of $X$. The proof of the proposition now follows from the fact that the rule $$(x,y) \mapsto (z,w)=(x^g,y)$$defines a regular covering map of degree $g$ from $C_{g,24}$ to the curve $w^3 = z^2-z.$

\begin{rema} The ``extra" involution of $C_{g,24}$ is given by the map $$\mathbf{s} : (x,y)\mapsto (1/x, -y/x^g).$$
\end{rema}

\section{Proof of Theorem \ref{gpar}}
Let $p \geqslant 7$ be a prime number, and let $X$ be a compact Riemann surface of genus $2p$ endowed with a non-cyclic group of automorphisms $G$ of order $6p$. If $P$ is the normal Sylow $p$-subgroup of $G$ then $\bar{G}=G/P$ is isomorphic to either $\mathbb{Z}_6$ or $\mathbf{D}_3.$ If we apply the Riemann-Hurwitz formula to the regular covering map $\pi_P: X \to Y:=X/P$ then $$4p-2 = 2p(h-1) +s(p-1),$$where $h$ is the genus of $Y$ and  $s$ is the number of fixed points of $X$ by $P$. A routine computation shows that, as $p \geqslant 7,$ $h=s=2$. We denote such fixed points by $x_1$ and $x_2$, and  their images by $\pi_P$ by $y_1$ and $y_2$. Observe that if $\bar{G}$ acts on $Y$ fixing $y_1$ and $y_2$, then $\bar{G} \cong \mathbb{Z}_6$ and the $G$-stabiliser of $x_i$ has order $6p$ for $i=1,2.$ Therefore, $G$ is cyclic. 

\s

So, we can assume that $y_1$ and $y_2$ are not fixed by the action of $\bar{G}.$ Then these points are permuted and the $\bar{G}$-stabilizer of $y_i$ has order three, for $i=1,2$. By arguing as above, we see that the signature of the action of $G$ on $X$ is either \begin{equation}\label{sgni}\sigma_1=(0; 6,6,3p) \,\mbox{ or } \,  \sigma_2=(0;2,2,3,3p).\end{equation}
 
Observe that, for both possible signatures, the group $G$ has a normal cyclic subgroup $N$ of order $3p,$ and therefore the Schur-Zassenhaus theorem implies that $$G \cong N \rtimes H \cong \mathbb{Z}_{3p} \rtimes \mathbb{Z}_2,$$where $H \cong G/N \cong \mathbb{Z}_2$ is a Sylow $2$-subgroup of $G.$ It follows that \begin{equation}\label{presedi}G \cong \langle a, b : a^{3p}=b^2=1, bab=a^r\rangle\end{equation}for some $r \in \mathbb{Z}$ that satisfies $r^2 \equiv 1 \mbox{ mod } 3p.$

\s

Assume that $\bar{G} \cong \mathbf{D}_3$. If write $u=a^p, v=a^3$ and $w=b$ then $$G \cong \langle u, v, w : u^3=v^p=w^2=1, [u,v]=1, wuw=u^{-1}, wvw=v^{t} \rangle$$where $t=\pm 1.$ We observe that necessarily $t=-1.$ In fact, if $t=1$ then $G \cong \mathbf{D}_3 \times \mathbb{Z}_p$ but this group does not have elements of order 6 (and hence cannot act with signature $\sigma_1$), and cannot be generated by elements of order 2 and 3 (and hence cannot act with signature $\sigma_2$). Thus, we assume $G \cong \mathbf{D}_{3p}.$ As this dihedral group does not have elements of order six, necessarily $G$ acts on $X$ with signature $\sigma_2.$ Conversely, if we consider the presentation \eqref{presedi} with $r=-1$ for $\mathbf{D}_{3p}$,  then the ske $(b, ba^{2p-1}, a^p, a)$ guarantees the existence of a  Riemann surface of genus $2p$ with a dihedral group of automorphisms of order $6p$ acting with signature $\sigma_2.$

\s

We now assume $\bar{G} \cong \mathbb{Z}_6.$ Observe that in such a case, $r$ in the presentation \eqref{presedi} satisfies $r \equiv 1 \mbox{ mod 3},$ and therefore $rp \equiv p \mbox{ mod } 3p.$ It follows that $a^p$ commutes with $K=\langle a^3, b\rangle.$ If $K$ is cyclic then $G \cong K \times \langle a^p\rangle$ is cyclic as well. However, an abelian group of order $6p$ cannot act with any of the signatures in \eqref{sgni}. We conclude that necessarily $K$ is dihedral and hence $G$ is isomorphic to $\mathbf{D}_{p} \times \mathbb{Z}_3.$ Conversely, we claim that both signatures $\sigma_1$ and $\sigma_2$ are realised by this group. Indeed, if we write $r:=a^3$, $s:=b$ and $t:=a^p$ then the presentation \eqref{presedi} turns into \begin{equation}\label{directo}G \cong \mathbf{D}_p \times \mathbb{Z}_3=\langle r,s,t: r^p=s^2=(sr)^2=t^3=[t,r]=[t,s]=1\rangle.\end{equation}The skes $(st,srt,r^{-1}t)$ and $(s,sr,t,(rt)^{-1})$ ensure the existence of Riemann surfaces of genus $2p$ with a group of automorphisms isomorphic to $\mathbf{D}_p \times \mathbb{Z}_3$ acting with signature $\sigma_1$ and $\sigma_2$ respectively. 

\s

Summing up, the compact Riemann surfaces $X$ of genus $2p$ endowed with a group of automorphisms $G$ of order $6p$ split into four cases.
\begin{enumerate}
\item $G \cong \mathbb{Z}_{6p}$ acting with signature $\sigma_0=(0; 3,6p,6p).$ 
\item $G \cong \mathbf{D}_{3p}$ acting with signature $\sigma_2=(0;2,2,3,3p).$ 
\item $G \cong\mathbf{D}_{p} \times \mathbb{Z}_3$ acting with signature $\sigma_2=(0;2,2,3,3p).$ 
\item $G \cong\mathbf{D}_{p} \times \mathbb{Z}_3$ acting with signature $\sigma_1=(0;6,6,3p).$ 
\end{enumerate}

\s

We denote by $X_1$ the (unique) Riemann surface corresponding to case (1).

\s

{\it Claim.} Case (4) gives rise to a unique Riemann surface, isomorphic to $X_1.$  

\s

To prove the first statement of the claim, notice that, in terms of the presentation \eqref{directo}, the elements of $G=\mathbf{D}_{p} \times \mathbb{Z}_3$ of order six are $sr^it^{\pm 1}$ where $i \in \{0, \ldots, p-1\},$ whereas the elements of order $3p$ are  $r^kt^{\pm 1}$ where $k \in \{1, \ldots, p-1\}.$ Thus, each ske representing an action of $G$ with signature $\sigma_1$ is of the form \begin{equation}\label{skess}(sr^{i+k}t^{\pm}, sr^it^{\pm}, r^kt^{\pm}).\end{equation}Observe that, up to an automorphism of the form $r \mapsto r^u$, we can assume $k=1$. In addition, if $i \neq 0$ then we apply the automorphism $s \mapsto sr^{-i}$ to assume that $i=0$. The two choices in $t^{\pm}$ are clearly equivalent by inversion. Hence, the skes \eqref{skess} are equivalent to $$\theta_0=(srt,st,rt),$$ whence the uniqueness follows. We denote this Riemann surface by $X_4$.

In order to prove that $X_4$ is isomorphic to $X_1,$  consider the group $$\tilde{G}=\mathbf{D}_p \times \mathbb{Z}_6=\langle R, S, T : R^p=S^2=T^6=(SR)^2=[R,T]=[S,T]=1\rangle.$$The ske of signature $\sigma_{14}=(0; 2,6,6p)$ given by $$\Theta_{14}=(ST^3, SR^{-1}T^2, RT)$$guarantees the existence of a Riemann surface $X_{14}$ of genus $2p$ endowed with an action of $\tilde{G}$ with signature $\sigma_{14}.$ 
We now consider the following normal subgroups of $\tilde{G}$ $$\tilde{H}_1=\langle RT\rangle \cong \mathbb{Z}_{6p} \mbox{ and } \tilde{H}_2=\langle S, R, T^2\rangle \cong \mathbf{D}_p \times \mathbb{Z}_3,$$and the restriction of the action $\Theta_{14}$ to these subgroups. The equalities $$\langle ST^3\rangle \cap \tilde{H}_1=\{1\} \mbox{ and } \langle SR^{-1}T^2\rangle \cap \tilde{H}_1=\langle T^2\rangle$$show that the signature of the induced action of $\tilde{H}_1$ on $X_{14}$ is $(0; 3,6p,6p)$. Consequently, the uniqueness of $X_1$ allows us to conclude that $X_{14} \cong X_1.$ In a very similar way, it can be proved that $X_{14} \cong X_4$, thanks to the afore proved uniqueness of $X_4$.
 
\s

The cases (2) and (3) correspond to complex one-dimensional families; we denote them by $\mathcal{F}_{g,1}$ and $\mathcal{F}_{g,2}$ respectively. We claim that these families are irreducible subvarieties of the moduli space $\mathscr{M}_{2p}.$ Indeed, by \cite{B90} and \cite{GG92}, it suffices to prove that the families are equisymmetric, namely, that there is a unique topological class of actions of $G$, for each family. For the former family, if we write $$G \cong \mathbf{D}_{3p}= \langle a, b : a^{3p}=b^2=1, bab=a^{-1}\rangle$$ then   each ske representing the action of $G$ on $X \in \mathcal{F}_{g,1}$ is equivalent to either  \begin{equation*}\label{naranjo}(b,ba^{2p-1}, a^{p}, a) \mbox{ or }(b,ba^{p-1}, a^{2p}, a).\end{equation*}However, these skes are equivalent under the braid transformation $$(x_1,x_2, x_3, x_4) \mapsto (x_1, x_3^{-1}x_2x_3, (x_2x_3)^{-1}x_3(x_2x_3), x_4)$$(see \cite[Proposition 2.5]{B90} for more details). Analogously,  in terms of the presentation \eqref{directo}, the action $\mathbf{D}_p \times \mathbb{Z}_3$ on each $X \in \mathcal{F}_{g,2}$ is represented by the ske $(s,sr^{-1}, t, rt^2)$. 

\s

The fact that the signature $(0; 2, 2,3,3p)$ is maximal (see \cite{Sing72}) implies that, up to finitely many exceptions, the automorphism group of each member of these families is isomorphic to $\mathbf{D}_{3p}$ and $\mathbf{D}_{p}\times \mathbb{Z}_3$ respectively. 
\s

We now proceed to determine the algebraic model for the members of $\mathcal{F}_{g,1}$.  First, observe that each $X \in \mathcal{F}_{g,1}$ is trigonal. In fact,  the regular covering map $\pi : X \to X/\langle a^p \rangle \cong \mathbb{P}^1$ ramifies over $2p+2$ values. We denote these points  by \begin{equation}\label{calef}\alpha_1, \ldots, \alpha_p, \beta_1, \ldots, \beta_p, \gamma_1, \gamma_2.\end{equation}It follows that $X$ is represented by the algebraic curve $$y^3=\Pi_{i=1}^p(x-\alpha_i)^{u_i}\Pi_{i=1}^p(x-\beta_i)^{v_i}(x-\gamma_1)^{\epsilon_1}(x-\gamma_2)^{\epsilon_2}$$where $u_i, v_i, \epsilon_i \in \{1,2\}$. Observe that  $a^3$ induces an automorphisms of $\mathbb{P}^1$ of order $p$ that keeps the set \eqref{calef} invariant. This symmetry allows us to assume that  $$\alpha_i=\omega_p^i, \,\beta_i=\lambda \omega_p^i, \,\gamma_1=0 \mbox{ and } \gamma_2=\infty,$$for some $\lambda \in \mathbb{C}$ such that $\lambda^p \neq 0,1.$ Consequently, $X$ is represented by the  curve $$y^3=x^{\epsilon}(x^p-1)^{u}(x^p-\eta)^{v} \mbox{ where } \eta:=\lambda^p.$$The exponents above can be determined from the rotation constants of $a^p$ at each one of its fixed points. Let $P=\pi^{-1}(\alpha_1)$ and $P'=\pi^{-1}(\beta_1).$ Without loss of generality, it can be assumed that $\mbox{rot}(P,a^p)=\omega_3$ (see, for instance, \cite[\S3.1]{B22}) and therefore $$\mbox{rot}(P',a^p)= \mbox{rot}(b(P),a^p)=\mbox{rot}(P,a^{-p})=\mbox{rot}(P,a^{p})^{-1}=\omega_3^2.$$It follows that $u=1$ and $v=2$. In a similar way, one sees that $\epsilon$ can be taken to be 1. Thus, we conclude that $X$ is algebraically represented by $$y^{3}=x(x^p-1)(x^p-\eta)^{2} \mbox{ where } \eta \neq 0, 1.$$

\s

An analog reasoning allows us to conclude that each $X \in \mathcal{F}_{g,2}$  is represented by  $$y^3=x^{\epsilon}(x^p-1)^{u}(x^p-\eta)^{v} \mbox{ where } \eta \neq 0, 1$$ for some suitably chosen $\epsilon, u, v \in \{1,2\}.$ In this case, the trigonal map is given by the action of $\langle t \rangle.$ As $t$ is central, we can assume that $u=v=1$ (equivalently, the rotation constant of $t$ at each one of its fixed points with $G$-stabiliser $\langle t \rangle$ can be assumed to be equal to $\omega_3$). In addition, if $P$ is the fixed point of $t$ with $G$-stabiliser $\langle rt^2\rangle$ corresponding to 0  then one can assume$$\mbox{rot}(P,rt^2)=\omega_{3p} \mbox{ and therefore } \mbox{rot}(P,t^{2p}) = \mbox{rot}(P,rt^2)^p=\omega_3.$$This shows that $\epsilon=1$ if $p \equiv 2\mbox{ mod }3$ and $\epsilon=2$ otherwise. Thus $X$ is isomorphic to $$y^3=x(x^p-1)(x^p-\eta) \, \mbox{ or }\,\, y^3=x^2(x^p-1)(x^p-\eta)$$respectively. Finally, by taking $\eta=-1$ above, we conclude that a point of $\mathcal{F}_{g,2}$ is represented by the curve $C_{g,2}$ if $p \equiv 2 \mbox{ mod 3}$ and by the curve $C_{g,4}$ otherwise.

\s

\begin{rema}
The behaviour for $p=2,3$ and 5 is, as expected, irregular. For instance, the alternating group $\mathbf{A}_4$ of order twelve acts on genus four with signature $(1; 2)$; there is a Riemann surface of genus ten with automorphism group of order thirty acting with signature $(0;5,6,30);$ and since $\mathbf{D}_3 \times \mathbb{Z}_3$ does not have elements of order nine, the family $\mathcal{F}_{6,2}$ does not exist.
\end{rema}

\section{Proof of Theorem \ref{final}}

After considering Theorems \ref{unicidad} and \ref{gpar}, the proof of the theorem follows directly after observing that the family $\mathcal{F}_{g,1}$ in the theorem above exists for each $p \geqslant 7$ non-necessarily prime. For even genera $g \leqslant 14$ the existence of a compact Riemann surface with a non-cyclic group of automorphisms of order $3g$ can be verified easily, for instance, by consulting the databases \cite{Conder}, \cite{Ka} and \cite{L}.

\section{Proof of Theorem \ref{gimpar}}
Let $p \geqslant 13$ be a prime number, and let $X$ be a compact Riemann surface of genus $4p-1$ endowed with a group of automorphisms $G$ of order $12p$. Assume that $G$ is non-cyclic. Let $P \cong \mathbb{Z}_p$ denote the normal Sylow $p$-subgroup of $G$. If we denote by $\gamma$ the genus of the quotient Riemann surface $Y=X/P$ then, by the Riemann-Hurwitz formula applied to the regular covering map $\pi: X \to Y$, one gets 
\begin{equation*}\label{rh4p}
8p-4 = 2p(\gamma-1) +s(p-1),
\end{equation*}where $s$ denotes the number of fixed points of the action of $P$. It follows that $\gamma \leqslant 4$. Furthermore, as $p\geqslant 13, $ it is not difficult to verify that $\gamma = 3$ and $s=4$. 

Since $Y$ has genus three and is equipped with a group of automorphisms $K:=G/P$ of order twelve, the list of possibilities for the algebraic structure of $K$ and its signature is known. We summarise this information below (see, for instance, \cite{B91}, \cite{Conder} and \cite{L}).
\s

  \begin{center}
\begin{tabular}{|c||c|c|c|c|c|}  
\hline
$K$ & $\mathbb{Z}_{12}$ & $\mathbb{Z}_{12}$ & $\mathbb{Z}_3 \rtimes \mathbb{Z}_4$& $\mathbf{D}_6$ &  $\mathbf{A}_4$\\  \hline
  signature &  $(0;2,12,12)$ & $(0;3,4,12)$ & $(0;4,4,6)$ & $(0;2,2,2,6)$ & $(0;2,2,3,3)$\\ 

\hline
\end{tabular}
\end{center}

\s

Let $\Omega \subset Y$ denote the set of branch values of $\pi.$ Note that $K$ must keep $\Omega$ set-wise invariant. The fact that $G$ is assumed to be non-cyclic implies that no element in $\Omega$ is fixed by $K$, and therefore $\Omega$ splits into two orbits of length two or the action of $K$ is transitive. Observe that the former case is impossible, since otherwise it would imply that 6 appears twice in the signature of the action of $K$ on $Y$. It then follows that the action of $K$ on $\Omega$ is transitive and therefore 3 appears in the signature of the action of $K$ on $Y$. Hence, one of the following statements holds.

\begin{enumerate}
\item $K \cong \mathbb{Z}_{12}$ acts on $Y$ with signature 
$(0;3,4,12)$, and therefore the signature of the action of $G$ on $X$ is $(0; 4,12,3p).$ 
\item $K \cong \mathbf{A}_{4}$ acts on $Y$ with signature $(0;2,2,3,3)$, and therefore the signature of the action of $G$ on $X$ is $(0;2,2,3,3p).$ 
\end{enumerate}

\s
We now proceed to consider these two cases separately.
\s

{\bf Case 1.} Assume that the action of $G$ on $X$ has signature $(0; 4,12,3p)$ and that the action is represented by the ske $\theta=(a,b,c).$ Observe that the Sylow $2$-subgroups of $G$ are cyclic, and that $\langle c \rangle \cong \mathbb{Z}_{3p}$ is a normal subgroup of $G$. It follows that $$G \cong \mathbb{Z}_{3p} \rtimes \mathbb{Z}_4=\langle A,B: A^{3p}=B^{4}=1, BAB^{-1}=A^r\rangle$$for some $r.$ Note that as $G$ has an element of order twelve, $[A^p,B]=1.$ If we write $u=A^p, x=A^3$ and $w=B$ then $$G \cong \langle u,x,w: u^3=x^p=w^4=1, [u,x]=1, wuw^{-1}=u^{t_1}, wxw^{-1}=x^{t_2}\rangle$$where $t_1=\pm 1$ and $t_2^4 \equiv 1 \mbox{ mod }p.$ As $[A^p,B]=1,$ we have that $t_1=1.$ In addition,  as $p \equiv -1 \mbox{ mod }4,$ one has that $t_2=\pm 1.$ Note that $t_2=-1$ since otherwise $G \cong \mathbb{Z}_{3p} \times \mathbb{Z}_4$ and this group does not act with signature $(0; 4,12,3p)$. After noticing that $y:=uw$ has order $12$ and $yxy^{-1}=x^{-1}$, we conclude that \begin{equation}\label{barao}G \cong \mathbb{Z}_{p} \rtimes_{2} \mathbb{Z}_{12}=\langle x,y: x^{p}=y^{12}=1, yxy^{-1}=x^{-1}\rangle.\end{equation}The elements of $G$ of order 4 of $G$ are $x^ny^{\pm 3},$ the ones of order 12  are $x^ny^{\pm 1}, x^ny^{\pm 5}$, and the ones of order $3p$ are $x^my^{\pm 4}$, where $n \in \{1, \ldots, p\}$ and $m \in \{1, \ldots, p-1\}.$ After considering the automorphism of $G$ given by $y \mapsto y^{-1},$ one sees that $a$ can be assumed to be of the form $x^ny^3$. Consequently, since $abc=1$, the ske $\theta$ has the form $$(x^ny^3,  x^{n'}y^5,x^{m}y^4) \, \mbox{ or } \, (x^ny^3,  x^{n'}y, x^{m}y^8)$$where $n,n' \in \{1, \ldots, p\}$ and $m \in \{1, \ldots, p-1\}.$ Observe that the automorphism of $G$ given by $y \mapsto y^5$ allows us to disregard the skes of the latter form. In addition, since $m \neq 0$, by considering an automorphism of $G$ of the form $x \mapsto x^u,$ we obtain that $m$ can be assumed to be 1. Hence, $\theta$ is equivalent to $$\theta_n:=(x^{n}y^3, x^{n+1}y^5, xy^4) \mbox{ for some } n \in \{1,\ldots, p\}.$$The automorphism of $G$ given by $y \mapsto xy$ shows that $\theta_n $ and $\theta_{n+1}$ are equivalent for each $n.$ It then follows that $$\theta \cong (x^{-1}y^3, y^5, xy^4)$$and therefore $X$ is uniquely determined by $\theta$, up to isomorphism. 

\s
Observe that the regular covering map $$\pi_H:X \to X/H \,\mbox{ where }\, H=\langle xy^4\rangle \cong \mathbb{Z}_{3p}$$ramifies over five values, one marked with 3 and remaining ones marked with $3p.$ Thus, $X/H$ is rational and  $X$ is $3p$-gonal. Since $X/H$ admits an action of $G/H \cong \mathbb{Z}_4$ in such a way that $(X/H)/(G/H) \cong X/G$, one sees that, up to a M\"{o}bius transformation, the branch values of $\pi_H$ marked with $3p$ are $1,i,-1,-i$ and the branch value marked with $3$ is $\infty.$ It follows that $X$ is isomorphic to the normalisation of the algebraic curve \begin{equation*}\label{lapiz}y^{3p}=(x-1)^{m_1}(x-i)^{m_2}(x+1)^{m_3}(x+i)^{m_4} \mbox{ for some } m_i \in \{1, \ldots, 3p-1\}.\end{equation*}

We claim that we can assume $m_1=m_3=1$ and $m_2=m_4=s$ where $s \in \{2, \ldots, 3p-1\}$ is the unique integer satisfying $$
  s \equiv -1 \mbox{ mod } p  \mbox{ and }
  s \equiv 1  \mbox{ mod } 3.$$
  
  In order to prove the claim first notice that the genus of the Riemann surface represented by the curve \begin{equation}\label{campinas}y^{3p}=(x^2-1)(x^2+1)^s\end{equation}is $4p-1.$ Observe that $s^2 \equiv 1 \mbox{ mod } 3p$ and therefore $s^2-1=3pL$ for some integer $L.$ If we define $m:=s-L$ then $3p(m-s)+s^2=1$ (note that $m$ is odd). We define $$\mathbf{a}(x,y)=(x, \omega_{3p}y) \mbox{ and } \mathbf{b}(x,y)=(ix,(-1)^{s+1}(x^2+1)^{m-s}y^s),$$ and a routine computation shows that $\mathbf{a}$ and $\mathbf{b}$ are automorphisms of \eqref{campinas}. If we write $\mathbf{x}:=\mathbf{a}^3$ and  $\mathbf{y}:=\mathbf{a}^p\mathbf{b}$ then $\langle \mathbf{x}, \mathbf{y}\rangle$ is isomorphic to \eqref{barao}. The conclusion follows from the uniqueness proved above.
  
\s

{\bf Case 2.} Assume that the action of $G$ on $X$ has signature $(0; 2,2,3,3p)$. Since $G$ has a normal subgroup $P$ of order $p$ and $G/P \cong {\bf A}_4$ we see that $G \cong \mathbb{Z}_p \rtimes {\bf A}_4.$ Observe that $G$ cannot be isomorphic to the direct product $\mathbb{Z}_p \times {\bf A}_4$ because this group is not generated by two involutions and one element of order three. If we write $$\mathbf{A}_4=\langle a,b,t : a^2=b^2=(ab)^2=t^3=1, tat^{-1}=b, tbt^{-1}=ab\rangle$$then $G=\langle z: z^p=1\rangle \rtimes {\bf A}_4$ where $$aza=z^{\pm 1}, bzb=z^{\pm 1} \mbox{ and }tzt^{-1}=z^{n} \, \mbox{ where } \, n^3 \equiv 1 \mbox{ mod } p.$$ Up to a permutation of the generators of $\langle a, b \rangle \cong \mathbb{Z}_2^2$ we can assume that either
$$aza=z, bzb=z^{-1} \,\, \mbox{ or } \,\, [z,a]=[z,b]=1.$$

The former case does not give rise to a group. In fact, since $tzt^{-1}=z^n$, the equality $tat^{-1}=b$ implies that $$z=tat^{-1}  z  (tat^{-1})^{-1}=bzb=z^{-1},$$a contradiction. It then follows that the latter case occurs and \begin{multline*}
G \cong \langle a,b,t,z : a^2=b^2=(ab)^2=t^3=z^p=1, tat^{-1}=b, tbt^{-1}=ab,\\ 
[a,z]=[b,z]=1, tzt^{-1}=z^n\rangle
\end{multline*}where $n$ is different from 1. However, this group does not have elements of order $3p.$

\begin{rema}
It would be worth to study if a result analogous to Theorem \ref{final} can be obtained for the case $3g+3.$ In contrast with the case $3g$, however, we conjecture that if g is odd and $g +1 \equiv 2 \mod 4$, the uniqueness property for Riemann surfaces with a group of automorphisms of order $3g+3$ should still hold (provided that such a surface exists). In fact, by using the same techniques as in Theorem \ref{gpar}, uniqueness can be proved for $g = 2p-1$, where $p\geqslant 7$ is a prime; we omit the details here. However, if $g = 2n-1$ for an odd non-prime integer $n$, the result, while extremely likely, is non-trivial to prove, and needs a deeper inspection. 
\end{rema}

\subsection*{Conflict of interest statement and Data availability statement} On behalf of all authors, the corresponding author states that there is no conflict of interest. This manuscript has no associated data.

%

\end{document}